
\input phyzzx
\input epsf.tex
\def\np{Nucl. Phys.}
\def\pl{Phys. Lett.}

\def\cmp{Comm. Math. Phys.}
\def\ijmp{Int. J. Mod. Phys.}
\def\jmp{J. Math. Phys.}

\def\inma{Invent. Math.}
\def\tam{Trans. Am. Math. Soc.}

\def\bams{Bull. AMS}
\def\am{Ann. of Math.}

\def\topo{Topology}

\def\knot{Journal of Knot Theory and Its Ramifications}
\def\half{{1\over 2}}

\def\ex{{\hbox{\rm e}}}

\def\tr{{\hbox{\rm Tr}}}

\def\vev{vacuum expectation value}

\def\erretres{{{\bf R}^3}}

\tolerance=500000
\overfullrule=0pt

\Pubnum={CTP-2434 \cr  US-FT-8-95 \cr hep-th/nnnmmyy}
\date={May, 1995}
\pubtype={}
\titlepage

\title{VASSILIEV INVARIANTS FOR TORUS KNOTS}
\author{M. Alvarez}
\address{Center for Theoretical Physics \break Massachusetts Institute of
Technology \break Cambridge, MA 02139, USA}
\andauthor{J.M.F. Labastida\foot{e-mail:
LABASTIDA@GAES.USC.ES}}
\address{Departamento de F\'\i sica de Part\'\i culas \break Universidade
de Santiago\break E-15706 Santiago de Compostela, Spain}

\abstract{Vassiliev invariants up to order six for arbitrary
torus knots $\{ n , m \}$, with $n$ and $m$ coprime integers, are computed.
These invariants are polynomials in $n$ and $m$ whose degree coincide with
their order. Furthermore, they turn out to be integer-valued in a normalization
previously proposed by the authors.}

\endpage
\pagenumber=1

\chapter{Introduction}

Vassiliev invariants \REF\vassi{V.A. Vassiliev, ``Cohomology of knot spaces"
 in {\it Theory of Singularities
and its applications}, (V.I. Arnold, ed.), Amer. Math.
Soc., Providence, RI, 1990, 23}
\REF\vassidos{V.A. Vassiliev, ``Topology of complements to discriminants
and loop spaces"  in {\it Theory of Singularities and its applications},
(V.I. Arnold, ed.), Amer. Math. Soc., Providence, RI, 1990, 23}
\REF\vassitres{V.A. Vassiliev, ``Complements of discriminants of smooth
maps: topology and applications", Translations of Mathematical Monographs,
vol. 98, AMS, 1992}
[\vassi,\vassidos,\vassitres]
seem to be a very promising set of knot
invariants to classify knot types.
Since the discovery
of their formulation in terms of inductive relations for singular
knots \REF\birlin{J.S. Birman and X.S. Lin\journal\inma&111(93)225}
\REF\birman{J.S. Birman\journal\bams&28(93)253}
 [\birlin,\birman], and of
their relation to knot invariants based on quantum groups or in
Chern-Simons gauge theory
\REF\drorcon{D. Bar-Natan, ``Weights of Feynman diagrams and the
Vassiliev knot invariants", preprint, 1991}
\REF\drortesis{D. Bar-Natan, ``Perturbative aspects of Chern-Simons
topological quantum field theory", Ph. D. Thesis, Princeton University,
1991}
\REF\lin{X.S. Lin, Vertex models, quantum groups and Vassiliev knot
invariants, Columbia University preprint, 1991}
[\drorcon,\drortesis,\birlin,\lin,\birman], several  works have been
performed to analyze Vassiliev invariants in both frameworks
\REF\konse{M. Kontsevich, {\sl Advances in Soviet Math.} {\bf 16}, Part 2
(1993) 137} \REF\drortopo{D. Bar-Natan\journal\topo&34(95)423}  \REF\piuni{S.
Piunikhin\journal\knot&4(95)163} \REF\numbers{M. Alvarez and J.M.F.
Labastida\journal\np&B433(95)555; Erratum\journal\np&B441(95)403} \REF\haya{N.
Hayashin, ``Graph Invariants of Vassiliev Type and Application to 4D Quantum
Gravity", UT-Komaba 94-8, March 1995, q-alg/9503010}
[\konse,\drortopo,\piuni,\numbers,\haya]. In [\konse,\drortopo] it was  shown
that Vassiliev invariants can be understood in terms of representations of
chord diagrams whithout isolated chords modulo the so called 4T relations
(weight systems), and that using semi-simple Lie algebras  weight systems can
be constructed. It was also shown in [\drortopo], using Kontsevitch's
representation for Vassiliev invariants
\REF\konsedos{M. Kontsevich, ``Graphs, homotopical algebra and
low-dimensional topology, preprint, 1992} {[\konsedos], that the space of
weight systems is the same as the space of Vassiliev invariants. In [\piuni]
it was shown that these representations are precisely the ones underlying
quantum-group or Chern-Simons invariants.

Recently, we made
the observation [\numbers] that the generalization of the integral or
geometrical knot invariant first proposed in \REF\gmm{E. Guadagnini, M.
Martellini and M. Mintchev\journal\pl&B227(89)111
\journal\pl&B228(89)489\journal\np&B330(90)575} [\gmm] and further
analyzed in
 [\drortesis], as well as the invariant itself are Vassiliev
invariants. In [\numbers] we proposed an organization of those
geometrical invariants and we described a procedure for their calculation
from known polynomial knot invariants. This procedure was applied to
obtain Vassiliev knot invariants up to order six for all prime knots up
to six crossings. These geometrical invariants have also been studied by
Bott and Taubes  \REF\bot{R. Bott and C. Taubes\journal\jmp&35(94)5247}
[\bot] using a different approach. In this paper we use the techniques
used in [\numbers] to compute all Vassiliev invariants up to order six
for arbitrary  torus knots.  Torus knots are labelled by two
coprime integers $n$ and $m$, $(n,m)=1$, such that the torus knots $\{n,m\}$,
$\{m,n\}$, $\{-n,-m\}$ and $\{-m,-n\}$ are the same knot, and $\{n,m\}$ and
$\{n,-m\}$ are mirror images of each other. The resulting invariants are
polynomials in
$n$ and $m$. This is consistent with the characterization found in
\REF\dean{J. Dean\journal\knot&3(94)7}
\REF\trapp{R. Trapp\journal\knot&3(94)391} [\dean,\trapp] for torus knots of
the form $n=2$ and $m=2p+1$. These polynomials turn out to be integer-valued
when
$(n,m)=1$ after they are normalized using the normalization proposed in
[\numbers].

The paper is organized as follows. In sect. 2 we give a brief description
of the framework proposed in [\numbers] for Vassiliev invariants. In
sect. 3 we collect known results on polynomials invariants for torus
knots associated to different groups and representations which will be
used in our computations. In sect. 4 we calculate the Vassiliev
invariants for torus knots up to order 6. Finally, in sect. 5 we analyze
the properties of these invariants. An appendix deals with some technical
details.

\chapter{Vassiliev Invariants from Chern-Simons Gauge Theory}

In this section we will recall a few facts on Chern-Simons gauge theory
and we will summarize the formulation for Vassiliev invariants based on
Chern-Simons perturbation theory proposed in [\numbers]. Let us consider
a gauge group $G$ and a connection $A$ on ${\bf R^3}$. The Chern-Simons
action is defined as:
$$
S_k(A)={k\over 4\pi}\int_{\bf R^3} \tr(A\wedge dA + {2\over 3} A\wedge
 A\wedge A),
\eqn\action
$$
where ``$\tr$" denotes the trace in the fundamental representation of
$G$. Given a knot $C$, \ie, an embedding of $S^1$ into ${\bf R^3}$, we
define the Wilson line associated to $C$ carrying a representation $R$
of $G$ as:
$$
W_C^R=\tr_R\big[ {\hbox{\rm P}} \exp \oint A \big],
\eqn\wilson
$$
where ``P" stands for path ordered and the trace is to be taken
in the representation $R$. The vacuum expectation value of $W_C^R$
is defined as the following ratio of functional integrals:
$$
\langle W_C^R \rangle = {1\over Z_k} \int [DA] W_C^R \ex^{iS_k(A)},
\eqn\vev
$$
where $Z_k$ is the partition function:
$$
Z_k=\int [DA] \ex^{iS_k(A)}.
\eqn\parfun
$$

The theory based on the action \action\ possesses a gauge symmetry which
has to be fixed. In addition, one has to take into account that the
theory must be regularized due to the presence of divergent integrals
when performing the perturbative expansion of \vev. Regarding these two
problems we will follow the approach taken in [\numbers]. Namely, we
will work in the Landau gauge and we will assume that there exist a
regularization such that one can ignore the shift in $k$,
$k\rightarrow k+g^\vee$ (being $g^\vee$ the dual coxeter number of $G$),
and Feynman diagrams which contain higher-loop contributions to two and
three-point functions (we refer the reader to [\numbers] for the details
concerning this issue).

There is one more problem emanating from perturbation quantum field
theory which must be considered. Often, products of operators
$A_\mu(x) A_\nu(y)$ must be considered at the same point $x=y$, where
they are divergent. This leads to an ambiguity which is solved guided by
the topological nature of the theory. In the process one needs to
introduce a framing attached to the knot which is characterized by an
integer $q$. It was shown in
\REF\pert{M. Alvarez and J.M.F. Labastida\journal\np&B395(93)198} [\pert]
that to work in the standard framing $(q=0)$ is equivalent to ignore
diagrams containing collapsible propagators in the sense explained
in [\pert,\numbers]. As in [\numbers] we will indeed work in the
standard framing.

As shown in [\numbers] the perturbative expansion of the vacuum
expectation value of the Wilson line \wilson\ has the form:
$$
\langle W_C^R \rangle = d(R) \sum_{i=0}^\infty \sum_{j=1}^{d_i}
\alpha_{ij} r_{ij} x^i,
\eqn\expansion
$$
where $x={2\pi i \over k}$, and $d(R)$ is the dimension of the
representation $R$. The factors $\alpha_{ij}$ and $r_{ij}$ in
\expansion\ incorporate all dependence dictated from the Feynman rules
apart from the dependence on $k$ which is contained in $x$. The power of
$x$, $i$, represents the order in perturbation theory. Of the two
factors, $r_{ij}$ and $\alpha_{ij}$, the first one contains all the
group-theoretical dependence, while the second all the geometrical
dependence. The quantity $d_i$ denotes  the number of
independent group structures $r_{ij}$ which appear at order $i$. The
first values of $d_i$, $\alpha_{ij}$ and $r_{ij}$ are: $\alpha_{0,1}=
r_{0,1}=1$, $d_0=1$ and $d_1=0$. Notice that we are dispensing with the
shift in $k$ ($x={2\pi i\over k}$), and therefore no diagrams with
higher-loop contributions to two and three-point functions should be
considered. In addition, there is no linear term in the expansion ($d_1=0$)
so that diagrams with collapsible propagators should be ignored in the
sense explained in [\numbers]. It was proven in [\numbers]
that the quantities $\alpha_{ij}$ are Vassiliev invariants of order $i$.

The group factors up to order $i=6$ have been computed in [\numbers] in
terms of the basic Casimir invariants. These group factors are
classified in two types, the ones which are not products of lower order
group factors:
$$
\eqalign{ r_{2,1}=&\sum_{l=1}^M C_3^{(l)}, \cr r_{3,1}=&\sum_{l=1}^M
\big(C_3^{(l)}\big)^2
\big(C_2^{(l)}\big)^{-1}, \cr r_{4,2}=&\sum_{l=1}^M \big(C_3^{(l)}\big)^3
\big(C_2^{(l)}\big)^{-2}, \cr r_{4,3}=&\sum_{l=1}^M C_4^{(l)}, \cr
r_{5,2}=&\sum_{l=1}^M
\big(C_3^{(l)}\big)^4 \big(C_2^{(l)}\big)^{-3}, \cr
r_{5,3}=&\sum_{l=1}^M C_4^{(l)} C_3^{(l)}
\big(C_2^{(l)}\big)^{-1}, \cr} \qquad\qquad
\eqalign{ r_{5,4}=&\sum_{l=1}^M C_5^{(l)}, \cr r_{6,5}=&\sum_{l=1}^M
\big(C_3^{(l)}\big)^5
\big(C_2^{(l)}\big)^{-4}, \cr r_{6,6}=&\sum_{l=1}^M C_4^{(l)}
\big(C_3^{(l)}\big)^2
\big(C_2^{(l)}\big)^{-2}, \cr r_{6,7}=&\sum_{l=1}^M C_5^{(l)} C_3^{(l)}
\big(C_2^{(l)}\big)^{-1},
\cr r_{6,8}=&\sum_{l=1}^M C_6^{1\,(l)}, \cr  r_{6,9}=&\sum_{l=1}^M
C_6^{2\,(l)}, \cr}
\eqn\factors
$$
and the ones which can be written as products of lower order ones:
$$
\eqalign{ r_{4,1}=& r_{2,1}^2, \cr
 r_{5,1}=& r_{2,1}r_{3,1}, \cr
 r_{6,1}=& r_{2,1}^3, \cr}
\qquad\qquad
\eqalign{
 r_{6,2}=& r_{3,1}^2, \cr
 r_{6,3}=& r_{2,1}r_{4,2}, \cr
 r_{6,4}=& r_{2,1}r_{4,3}. \cr}
\eqn\losotros
$$

\midinsert
\centerline{\epsfbox{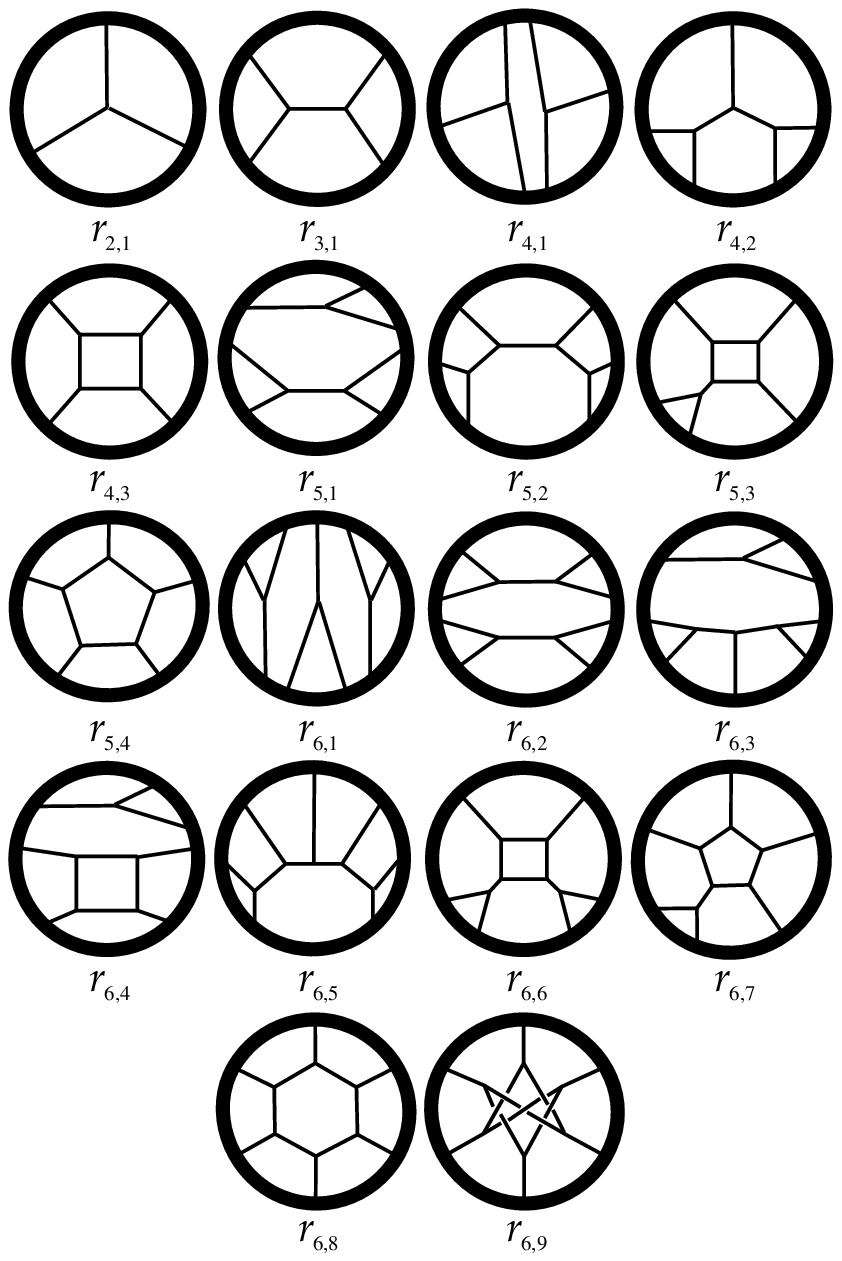}}
\centerline{Fig. 1. Diagrams associated to the independent group
factors.}
\endinsert

We will refer to the group factors in \factors\ as primitive group
factors, and we will denote by $\tilde d_i$ the number of primitive
group factors at order $i$. From \factors\ and \losotros\ follow that
the values of $d_i$ and $\tilde d_i$ for $i=1$ to $6$ are:
$$
\eqalign{
d =& 0, \,\,\,\, 1, \,\,\,\, 1, \,\,\,\, 3, \,\,\,\, 4, \,\,\,\, 9,\cr
\tilde d =& 0, \,\,\,\, 1, \,\,\,\, 1, \,\,\,\, 2, \,\,\,\, 3, \,\,\,\,
5.\cr}
\eqn\lasdes
$$
These dimensions are known up to $i=9$ [\drortopo].
The notation used in \factors\ is the following. We will assume that the
semi-simple group $G$ has the form $G=\otimes_{l=1}^n G^{(l)}$ where
$G^{(l)}$, $l=1,\dots,M$, are simple ones. The quantities
$C^{(l)}_2$, $C^{(l)}_3$, $C^{(l)}_4$, $C^{(l)}_5$, $C^{1(l)}_6$ and
$C^{2(l)}_6$ are traces of Casimir operators which for any simple group
$H$ have the
form (we remove the superindex $(l)$):
$$
\eqalign{
C_2 d(R) & = \tr(T_a T_a), \cr
C_3 d(R) & = - f_{abc} \tr(T_a T_b T_c), \cr
C_4 d(R)
&=f_{apq}f_{bqr}f_{crs}f_{dsp}\tr (T_a T_b T_c T_d), \cr C_5 d(R)
&=f_{apq}f_{bqr}f_{crs}f_{dst}f_{etp}\tr (T_a T_b T_c T_d T_e), \cr
C_6^1 d(R) &=f_{apq}f_{bqr}f_{crs}f_{dst}f_{etu}f_{rup}\tr (T_a T_b T_c
T_d T_e T_r), \cr C_6^2 d(R)
&=f_{apq}f_{brs}f_{ctp}f_{dur}f_{eqs}f_{gtu}\tr (T_a T_b T_c T_d T_e
T_g), \cr}
\eqn\casimires
$$
where $T_a$, $a=1,\dots,{\hbox{\rm dim}}(H)$, are the generators of the
simple group $H$ in the representation $R$, and $f_{abc}$
are the structure constants. The diagrams associated to these
Casimirs are presented in Fig. 2 while the diagrams associated to the
groups factors \factors\ and \losotros\ are contained in Fig. 1.

\midinsert
\centerline{\epsfbox{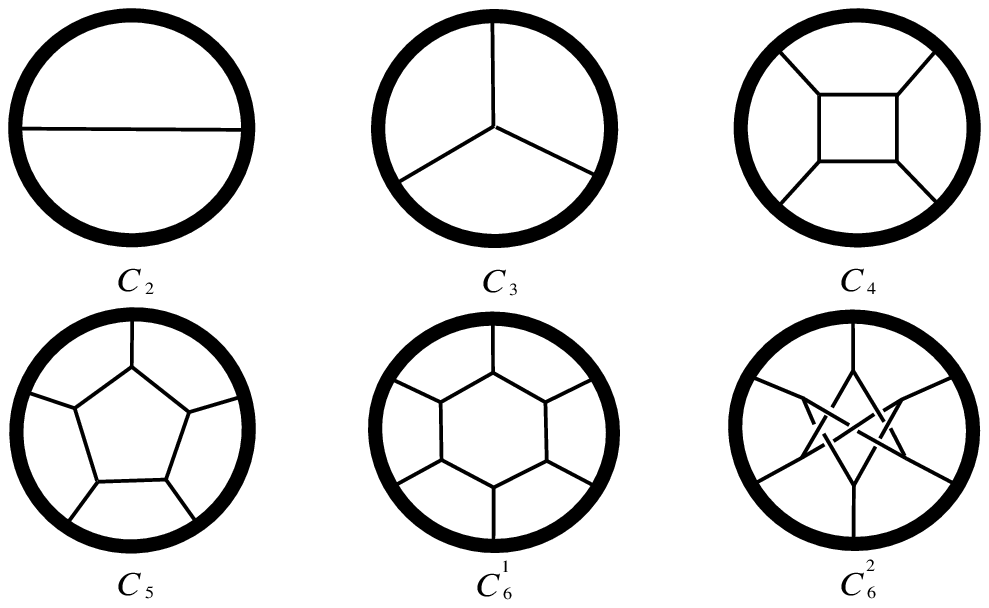}}
\centerline{Fig. 2. Diagramas associated to Casimir operators.}
\endinsert

As in [\numbers] we will use information from knot invariants for $SU(N)$
and $SO(N)$ in the fundamental representation, and for
$SU(2)$ in an arbitrary representation of spin $j/2$. The form of the
Casimirs \casimires\ for these cases, taken in a normalization such that
for the fundamental representation,
$$
\tr(T_a T_b) = - {1\over 2} \delta_{ab},
\eqn\norma
$$
is [\numbers]:

\noindent$\bullet$ $SU(N)_f$:
$$
\eqalign{C_2=& -{1\over 2N}(N^2-1),\cr
                C_3=& -{1\over 4}(N^2-1),\cr
                C_4=& {1\over 16}(N^2-1)(N^2+2),\cr
                C_5=& {1\over 32}N(N^2-1)(N^2+1),\cr
                C_6^1=& {1\over 64}(N^2-1)(N^4+N^2+2),\cr
                C_6^2=& {1\over 64}(N^2-1)(3N^2-2);\cr}
                \eqn\casivalsun
$$
\noindent$\bullet$ $SO(N)_f$:
$$\eqalign{C_2=& -{1\over 4}(N-1),\cr
                C_3=& -{1\over 16}(N-1)(N-2),\cr
                C_4=& {1\over 256}(N-1)(N-2)(N^2-5N+10),\cr
                C_5=& {1\over 1024}(N-1)(N-2)(N^3-7N^2+17N-10),\cr
                C_6^1=& {1\over 4096}(N-1)(N-2)(N^2-7N+14)(N^2-2N+3),\cr
                C_6^2=& {1\over 4096}(N-1)(N-2)(N-3)(7N-18);\cr}
                \eqn\casivalson
$$
\noindent$\bullet$ $SU(2)_j$:
$$\eqalign{C_2=& -j(j+1),\cr
                C_3=& -j(j+1),\cr
                C_4=& 2j^2(j+1)^2,\cr
                C_5=& 3j^2(j+1)^2-j(j+1),\cr
                C_6^1=& 2j^3(j+1)^3+3j^2(j+1)^2-2j(j+1), \cr
                C_6^2=& -2j^3(j+1)^3+5j^2(j+1)^2-2j(j+1).\cr}
\eqn\valores
$$

Equations \factors, \losotros\ and \casimires\ contain
 all the information appearing on the right hand side of \expansion\ up
to order six except the Vassiliev invariants themselves
(the quantities $\alpha_{ij}$) which we intend to
compute. To solve for them we need to know the left hand side of
\expansion. This information is compiled in next section.

\chapter{HOMFLY, Kauffman, Jones, and Akutsu-Wadati Polynomials for Torus
Knots}

In this section we convenientely present the form of the polynomial
invariants for torus knots corresponding to the groups $SU(N)$ and
$SO(N)$ in the fundamental representation, and to the group $SU(2)$ in an
arbitrary representation of spin $j/2$. These quantities have been worked
out from different points of view in the last few years.

The polynomial invariant for $SU(N)$ in the fundamental representation
is the HOMFLY polynomial
\REF\jones{V.F.R. Jones \journal\bams&12(85)103}
\REF\jonesAM{V.F.R. Jones\journal\am &  126  (87) 335}
\REF\homflyp{P. Freyd, D. Yetter, J. Hoste, W.B.R. Lickorish, K. Millet
and A. Ocneanu \journal\bams&12(85)239}
[\jones,\homflyp,\jonesAM]. It was first computed for torus knots
in [\jonesAM], reobtained from quantum groups in
\REF\rosso{M. Rosso and V. Jones\journal\knot&2(93)97} [\rosso], and from
Chern-Simons gauge theory in
\REF\homlinks{J.M.F. Labastida and M. Mari\~no\journal\ijmp&A10(95)1045}
[\homlinks]. It has the form:
$$
\eqalign{
\tilde P_{n,m}((\lambda t)^{1\over 2},t^{1\over 2}-t^{-{1\over 2}}) & =
\Big ( {1-t \over 1-t^n }\Big) {\lambda^{{1\over2}(m-1)(n-1)} \over
\lambda t-1} \cr & \,\,\,
\times\sum_{p+i+1=n \atop p, i \ge 0} (-1)^i t^{mi + {1 \over 2} p(p+1)}
{\prod_{j=-p}^i (\lambda t-t^j) \over (i)! (p)!} \cr}
\eqn\nuria
$$
where $(x)=t^x-1$, $\lambda= t^{N-1}$, and $t=\exp(2\pi i/(k+g^\vee))$ with
$g^\vee=N$.

The polynomial invariant for $SO(N)$ in the fundamental representation
corresponds to the Kauffman polynomial
\REF\kauf{L.H. Kauffman\journal\tam&318(90)417} [\kauf]. It was first
computed for torus knots in
\REF\yokota{Y. Yokota\journal\topo&32(93)309} [\yokota] and reobtained
from Chern-Simons gauge theory in
\REF\esther{J.M.F. Labastida and E. P\'erez, ``A relation between the
Kauffman and the HOMFLY polynomials for torus knots", Santiago preprint,
May 1995} [\esther]. It has the following form:
$$
\eqalign{
\tilde F_{n,m}(\lambda,t^{1\over 2}-t^{-{1\over 2}}) = &{{
{[1]}\,\lambda^{nm}} \over { {[1] + [ 0;1 ] }  }} \times \Biggl(
 \sum\limits_{\gamma +\beta +1=n\atop \gamma,\beta\geq 0} t^{- {m \over 2}
(\beta - \gamma)} \lambda^{-m}
   (-1)^{\gamma}
 \times \bigg(  {1 \over { [n]}} + { 1 \over { [ \beta - \gamma ; 1 ]}}
\bigg) \cr
&\times {1 \over {[\beta]!\,\ [\gamma]!}} \times
 \prod\limits_{j=-\gamma }^{\beta } [j;1] + \cases{ 0 &n odd \cr 1 &n
even \cr}
 \Biggr), \cr}
\eqn\naxos
$$
where $[p]=t^{p\over 2}-t^{-{p\over 2}}$, $[p;q]=t^{p\over 2}\lambda^q
-t^{-{p\over 2}}\lambda^{-q}$,
$\lambda=t^{N-1\over 2}$, and $t=\exp(\pi i/(k+g^\vee))$ with $g^\vee=N-2$.

The polynomial invariant for $SU(2)$ in an arbitrary representation of
spin $j/2$ is known as the Jones polynomial [\jones] for
$j={1}$, and the
Akutsu-Wadati polynomial for $j>{1}$. Its form for torus knots
was first obtained in \REF\poli{J.M. Isidro, J.M.F. Labastida and A.V. Ramallo
\journal\np&B398(93)187} [\poli]. The framework utilized was based on
the formalism proposed in
\REF\lr{J.M.F. Labastida and A.V. Ramallo \journal\pl&B227(89)92
\journal\pl&B228(89)214} [\lr]. Later it was reobtained in
\REF\kguno{R.K. Kaul, T.R.
Govindarajan\journal\np&B380(92)293\journal\np&B393(93)392
\journal\np&B402(93)548} [\kguno]
also in the framework of Chern-Simons gauge theory but using a
different formalism. This invariant has the form:
$$
\tilde I_{n,m}^j(t)={t^{{j\over 2}(n-1)(m-1) }\over t^{j+1}-1}
\sum_{l=0}^j t^{n(1+ml)(j-l)}(t^{1+sl}-t^{m(j-l)}),
\eqn\tdos
$$
where $t=\exp(2 \pi i/(k+2))$.

The invariants listed in \nuria, \naxos\ and \tdos\ contain enough
information to obtain the quantities $\alpha_{ij}$ in \expansion\ up to
order six. Actually, we must take into account that in \nuria, \naxos\
and \tdos\ the invariants have been normalized in such a way that for
the unknot they take the value one. In the expression \expansion\ we
need however the unnormalized invariants, \ie, the expressions for the
Wilson lines as obtained from \vev. These are easily obtained
multiplying by the value of the Wilson lines corresponding to the unknot
as dictated from Chern-Simons gauge theory. These expressions can be
extracted, for example, from
\REF\witCS{E. Witten\journal\cmp&121(89)351} [\witCS]
\REF\horne{J.H. Horne \journal\np&B334(90)669} and [\horne]
 for the cases of $SU(N)$ and $SO(N)$ in the fundamental
representation respectively, and from [\poli] for $SU(2)$ in a
representation of spin $j/2$. Taking into account these
facts it turns out that for the left-hand side of \expansion\ we must
consider: $$
P_{n,m}((\lambda t)^{1\over 2},t^{1\over 2}-t^{-{1\over 2}})  =
{(\lambda t)^{1\over 2} - (\lambda t)^{-{1\over 2}} \over
t^{1\over 2}-t^{-{1\over 2}}}
\tilde P_{n,m}((\lambda t)^{1\over 2},t^{1\over 2}-t^{-{1\over 2}}),
\eqn\nuriados
$$
where $(x)=t^x-1$, $\lambda= t^{N-1}$ and $t=\exp(2\pi i/(k+N))$,
$$
 F_{n,m}(\lambda,t^{1\over 2}-t^{-{1\over 2}})=
(1+{\lambda - \lambda^{-1} \over t^{1\over 2}-t^{-{1\over 2}}})
\tilde F_{n,m}(\lambda,t^{1\over 2}-t^{-{1\over 2}}),
\eqn\naxosdos
$$
where $\lambda=t^{N-1\over 2}$ and $t=\exp(\pi i/(k+N-2))$, and
$$
I_{n,m}^j(t) = {t^{j+1\over 2}-t^{-{j+1\over 2}}
\over t^{1\over 2}-t^{-{1\over 2}} }\tilde I_{n,m}^j(t),
\eqn\tdosdos
$$
where $t=\exp(2 \pi i/(k+2))$.

To obtain the linear equations satisfied by the $\alpha_{ij}$
one must appropriately expand the invariants \nuriados, \naxosdos\ and
\tdosdos. According to the form of $t$ entering in each of them, and
taking into account that one is dispensing with the shift, one has
to introduce the variable $x$ through:
$$
\eqalign{
t=&\ex^x, \,\,\,\,\,\,\,\,\,\, {\hbox{\rm for}} \,\,\,\, SU(N), \cr
t=&\ex^{x\over 2}, \,\,\,\,\,\,\,\,\,\, {\hbox{\rm for}} \,\,\,\,
SO(N).\cr}
\eqn\lastes
$$

\chapter{Vassiliev Invariants for Torus Knots up to order six}

In this section we obtain the expression of the Vassiliev invariants
${\alpha_{ij}}$ in \expansion\ up to order six. To carry this out one
must solve the equations resulting after considering the power series
expansion of the left hand side of \expansion\ once the variable $x$ has
been introduced as dictated by $t$ in \lastes.

In expanding \nuriados\ and \naxosdos\ one finds the slight complication
that $n$ appears in the limits of the sum and the product entering those
expressions. A similar complication shows up in \tdosdos\ due to the
presence of a sum whose upper limit is $j$.
To solve these difficulties we will proceed in the following
way. First notice that the invariants $P_{n,m}$, $F_{n,m}$ and $I_{n,m}$
are symmetric in $n$ and $m$. This implies that the coefficient of the
Taylor series expansion in $x$ must also be symmetric in $n$ and $m$.
On the other hand, for $n$ fixed in \nuriados\ and \naxosdos, and for $j$
fixed in \tdosdos, it is clear that the coefficients are polynomials in
$m$, and in $n$ and $m$ respectively. Actually, it was shown in [\dean,\trapp]
that if $n=2$ and $m=2p+1$ the polynomial in $p$ resulting as the
coefficient of $x^i$ is at most a polynomial of degree $i$. We can
therefore assume that the coefficient of $x^i$ in the Taylor series
expansion of \nuriados, \naxos\ and \tdosdos\ is a symmetric polynomial
in $n$ and $m$ of degrees $i$ in both $n$ and $m$.

One can actually state one more property  of the polynomials in $n$
and $m$ entering the coefficients of the Taylor series expansion of
\nuriados, \naxosdos\ and \tdosdos.
Under an inversion of space, a torus knot $\{n,m\}$ becomes the torus knot
$n,-m$. On the other hand, it is known that HOMFLY, Kauffman and
Akutsu-Wadati invariants which correspond to knots which are mirror
images of each other are related by the transformation $t \rightarrow
t^{-{1}}$. This implies that the coefficients of
the Taylor series expansions of
\nuriados, \naxosdos\ and \tdosdos\ are invariant under the change
$m\rightarrow -m$ (or $n\rightarrow -n$) and $x\rightarrow -x$. The
consequences of this symmetry are the following: the polynomial of the
coefficient of $x^i$ with $i$ even contains only even powers of $n$ and
$m$; the polynomial of the coefficient of $x^i$ with $i$ odd is of the
form $nm$ times a polynomial containing only even powers of $n$ and $m$.

The properties discussed so far are satisfied by both the
 unnormalized  invariants \nuriados, \naxosdos\ and \tdosdos,
and the normalized ones \nuria, \naxos\ and \tdos.
There is, however, one important property shared by the normalized invariants
which do not have the unnormalized ones and that suggests that one should
analyze first those. One can verify from \nuria, \naxos\ and \tdos\ that
$P_{1,\pm m}=F_{1,\pm m}=I^j_{1,\pm m}=1$. This implies that
$\pm 1$ must be roots of all the polynomials entering the coefficients
of $x^i$ for $i>0$ in both, $n$ and $m$, in the corresponding power
series. To exploit this property, instead
of considering the expansion \expansion\ for the unormalized  Wilson
lines we will consider the corresponding expansion for the normalized
ones,
$$
\langle \tilde W_C^R \rangle = {\langle W_C^R \rangle \over
\langle W_0^R \rangle} = \sum_{i=0}^\infty \sum_{j=1}^{d_i}
\tilde\alpha_{ij} r_{ij} x^i,
\eqn\expansiondos
$$
where $\langle W_0^R \rangle$ denotes the vacuum expectation
value for the unknot in the representation $R$.
Notice that we are using the same notation as in [\numbers]. We will
compute the quantities $\tilde \alpha_{ij}$ up to order six. Clearly,
from that result plus the fact the for the unknot [\numbers]:
$$
\eqalign{
&\langle W_0^R \rangle = \,d(R)\, \Big({{r_{1,2}}\over 6}\,{x^2}
+ ( {{{{r_{1,2}}^2}}\over {72}} + {{r_{2,4}}\over {360}} -
       {{r_{3,4}}\over {360}} ) \,{x^4}
+ ( {-{{{r_{1,2}}^3}}\over {1296}} -
{{r_{1,2}\,r_{2,4}}\over {2160}} \cr &
\,\,\,\,\,\,\,\,\,\,\,\,\,\,\,\,\,\,\,\,
   + {{r_{1,2}\,r_{3,4}}\over {2160}} - {{r_{5,6}}\over {15120}} +
     {{r_{6,6}}\over {3780}} - {{r_{7,6}}\over {11340}} +
{{r_{8,6}}\over {9072}} -
       {{r_{9,6}}\over {15120}} ) \, {x^6} + O(x^8)\Big), \cr}
\eqn\unknot
$$
one easily obtains the expression for the $\alpha_{ij}$.

{}From the properties described in the previous paragraphs follow that if
$i$ is even, the polynomial corresponding to the coefficient of $x^i$ in
the Taylor series expansion of \nuria, \naxos\ and \tdos\ has
the form $(n^2-1)(m^2-1)$ times a polynomial in $m^2$ and  $n^2$ which is
symmetric under $m^2 \leftrightarrow n^2$, and possesses at most degree
$i-2$ in both $n$ and $m$. For $i$ odd the corresponding polynomial has
the form $nm(n^2-1)(m^2-1)$ times a polynomial in $m^2$ and $n^2$ which is
symmetric under $m^2 \leftrightarrow n^2$, and possesses at most degree
$i-3$ in both $n$ and $m$.

Up to order $i=6$ the general form of the Taylor series expansions
of \nuria, \naxos\ and \tdos\ becomes:
$$
\eqalign{
 & 1 + (n^2-1)\,(m^2-1) \,{x^2}\,g_{2,1}  +
 n\,m\,(n^2-1)\,(m^2-1) \,{x^3}\,g_{3,1} \cr & +
  (n^2-1)\,(m^2-1) \,{x^4}\,
   ( g_{4,1} + ( {n^2} + {m^2} ) \,g_{4,2} + {n^2}\,{m^2}\,g_{4,3}
      ) \cr & +n\,m\,(n^2-1)\,(m^2-1) \,
   {x^5}\,( g_{5,1} + ( {n^2} + {m^2} ) \,g_{5,2} +
     {n^2}\,{m^2}\,g_{5,3} )  \cr &+
  (n^2-1)\,(m^2-1) \,{x^6}\,
   ( g_{6,1} + ( {n^2} + {m^2} ) \,g_{6,2} + {n^2}\,{m^2}\,g_{6,3} \cr &
\,\,\,\,\,\,\,\,\,\,\,\,\,\,\,\,\,\,\,\,\,\,
+   ( {n^4}\,{m^2} + {n^2}\,{m^4} ) \,g_{6,4} +
     ( {n^4} + {m^4} ) \,g_{6,5} + {n^4}\,{m^4}\,g_{6,6} ), \cr}
\eqn\ansatz
$$
where $g_{ij}$ are functions of $N$ for the cases of the HOMFLY and the
Kauffman invariants, and functions of the spin $j/2$ for the case of
Akutsu-Wadati.
The simplest way to determine these functions is to generate linear
relations among them computing the Taylor series expansions of the
left-hand side of \expansiondos\ for specific values of $n$ and $m$ up to
order six. In doing this one finds the following results:

\noindent$\bullet$ $SU(N)_f$:
$$
\eqalign{
g_{2,1}=& -{1\over 24}(N^2-1), \cr
g_{3,1}=& -{1\over 144}N(N^2-1), \cr
g_{4,1}=& {1\over 5760}(7N^4-10N^2+3), \cr
g_{4,2}=& -{1\over 5760}(3N^4-10N^2+7), \cr}
\qquad\qquad
\eqalign{
g_{4,3}=& -{1\over 1920}(N^4-1), \cr
g_{5,1}=& {1\over 86400}N(29N^4-60N^2+31), \cr
g_{5,2}=& -{1\over 86400}N(11N^4-40N^2+29), \cr
g_{5,3}=& -{1\over 86400}N(N^4-10N+11), \cr}
$$
$$
\eqalign{
g_{6,1}=& -{1\over 967680}(31N^6-49N^4+21N^2-3), \cr
g_{6,2}=& {1\over 483840}(9N^6-35N^4+35N^2-9), \cr
g_{6,3}=& {1\over 1451520}(55N^6-98N^4-35N^2+78), \cr
g_{6,4}=& -{1\over 1451520}(22N^6-77N^4+28N^2+27), \cr
g_{6,5}=& -{1\over 967680}(3N^6-21N^4+49N^2-31), \cr
g_{6,6}=& {1\over 2903040}(5N^6-49N^4+35N^2+9); \cr}
\eqn\piliuno
$$
\noindent$\bullet$ $SO(N)_f$:
$$
\eqalign{
g_{2,1}=& -{1\over 96}(N^2-3N+2), \cr
g_{3,1}=& -{1\over 1152}(N-2)^2(N-1), \cr
g_{4,1}=& {1\over 92160}(7N^4-45N^3+110N^2-120N+48), \cr
g_{4,2}=& -{1\over 92160}(3N^4-15N^3+20N^2-8), \cr}
$$
$$
\eqalign{
g_{4,3}=& -{1\over 92160}(3N^4-25N^3+70N^2-80N+32), \cr
g_{5,1}=& {1\over 5529600}(58N^5-469N^4+1455N^3-2120N^2+1412N-336), \cr
g_{5,2}=& -{1\over 5529600}(22N^5-141N^4+295N^3-180N^2-92N+96), \cr
g_{5,3}=& -{1\over 5529600}(2N^5-51N^4+245N^3-480N^2+428N-144), \cr}
$$
$$
\eqalign{
g_{6,1}=& -{1\over
61931520}(31N^6-315N^5+1358N^4-3150N^3+4116N^2\cr &
{\hbox{\hskip7.2cm}}-2856N+816), \cr
g_{6,2}=&
{1\over 61931520}(18N^6-147N^5+455N^4-630N^3+280N^2+168N-144), \cr
g_{6,3}=& {1\over
928972800}(550N^6-5768N^5+23443N^4-46865N^3+47740N^2 \cr &
{\hbox{\hskip7.2cm}}-22652N+3552),
\cr g_{6,4}=& -{1\over
928972800}(220N^6-1757N^5+5152N^4-6335N^3+1540N^2 \cr &
{\hbox{\hskip7.2cm}}+3052N-1872), \cr
g_{6,5}=& -{1\over 61931520}(3N^6-21N^5+42N^4-56N^2+32), \cr
g_{6,6}=& {1\over
928972800}(25N^6+112N^5-1442N^4+4585N^3-6860N^2\cr &
{\hbox{\hskip7.2cm}}+5068N-1488); \cr}
\eqn\pilidos
$$
\noindent$\bullet$ $SU(2)_j$:
$$
\eqalign{
g_{2,1}=& {1\over 6}A, \cr
g_{3,1}=& {1\over 18}A, \cr
g_{4,1}=& {1\over 360}A(7A-1), \cr
g_{4,2}=& -{1\over 360}A(3A+1), \cr}
\qquad\qquad
\eqalign{
g_{4,3}=& {1\over 360}A(7A+9), \cr
g_{5,1}=& {1\over 1080}A(10A-1), \cr
g_{5,2}=& -{1\over 1080}A(6A+3), \cr
g_{5,3}=& {1\over 1080}A(18A+15), \cr}
$$
$$
\eqalign{
g_{6,1}=& {1\over 75600}A(155A^2-55A^2+5), \cr
g_{6,2}=& -{1\over 75600}A(90A^2+20A-5), \cr
g_{6,3}=& {1\over 75600}A(260A^2+358A-9)), \cr
g_{6,4}=& -{1\over 75600}A(90A^2+342A+184), \cr
g_{6,5}=& {1\over 75600}A(15A^2+15A+5), \cr
g_{6,6}=& {1\over 75600}A(155A^2+1023A+691), \cr}
\eqn\pilitres
$$
where,
$$
A=-{1\over 4}j(j+2).
\eqn\laa
$$

We have all the ingredients to compute the Vassiliev invariants for torus
knots up to order six. Equations \piliuno, \pilidos\ and \pilitres\
plugged in \ansatz\ constitute the left-hand side of \expansiondos. Since all
but the $\tilde\alpha_{ij}$ is known on the right-hand side of that equation we
can obtain linear relations for the $\tilde\alpha_{ij}$ considering
different groups and representations. As in [\numbers] we will take into
account the linear relations appearing from the consideration of $SU(N)$
and $SO(N)$ in their fundamental representations,  $SU(2)$ in a
representation of arbitrary spin $j/2$, and $SU(N)\times SU(2)$ in a
representation which is the product of the fundamental of $SU(N)$ times
the one of spin $j/2$ of $SU(2)$. For this last case we will use the
following property.  Denoting by
$\langle W^{R_l}_C\rangle$ and  $\langle W^R_C \rangle$ the vacuum
expectation values of Wilson lines based on the simple groups $G_l$ where
$G=\otimes_{l=1}^M G_l$, being the representation $R$ of $G$ a direct product
of representations $R_l$ of $G_l$, one has,
$$
\langle W^R_C \rangle =\prod_{l=1}^M \langle W^{R_l}_C \rangle.
\eqn\factorization
$$
This follows directly from the factorization of both the partition
function and the Wilson line operator [\numbers].

Proceeding in the way described above one finds 5
equations for $\tilde\alpha_{2,1}$, 5 equations for $\tilde\alpha_{3,1}$,
12 equations for $\tilde\alpha_{4,1}$, $\tilde\alpha_{4,2}$ and
$\tilde\alpha_{4,3}$, 15 equations for
$\tilde\alpha_{5,1},\cdots,\tilde\alpha_{5,4}$, and 20 equations for
$\tilde\alpha_{6,1},\cdots,\tilde\alpha_{6,9}$. These equations have a unique
solution which takes the form:
$$
\eqalign{
\tilde\alpha_{2,1}=& {1 \over 6}  \,( {n^2}-1 ) \,({m^2}-1),  \cr
\tilde\alpha_{3,1}=& {1 \over 18} \,n\,m\,(  {n^2}-1 ) \,(  {m^2}-1 ), \cr
\tilde\alpha_{4,1}=& {1 \over 72} \,(n^2-1)^2\,(m^2-1)^2, \cr
\tilde\alpha_{4,2}=& {1 \over 360} \,( {m^2}-1 ) \, (n^2-1) \,
     ( 9\,{n^4}\,{m^2}-m^2-n^2-1 ), \cr
\tilde\alpha_{4,3}=& {1 \over 360} \,(  {n^4}-1 ) \,(  {m^4}-1 ), \cr}
$$
$$
\eqalign{
\tilde\alpha_{5,1}=& {1 \over 108}\,n \,m \,(n^2-1)^2\,(m^2-1)^2, \cr
\tilde\alpha_{5,2}=& {1 \over 5400} \,n\,m \,( n^2-1 )\, ( m^2-1) \,
     ( 69\,{n^2}\,{m^2}  - 21\,\big({n^2} + {m^2})-11\big), \cr
\tilde\alpha_{5,3}=& {1 \over 5400} \,n\,m \,( n^2-1 )\, ( m^2-1) \,
     ( 11\,{n^2}\,{m^2} + {n^2} + {m^2}-9), \cr
\tilde\alpha_{5,4}=& {1 \over 900} \,n\,m\,(  {n^4}-1 ) \,( {m^4}-1 ), \cr}
$$
$$
\eqalign{
\tilde\alpha_{6,1}=& {1 \over 1296} \,(n^2-1)^3\,(m^2-1)^3, \cr
\tilde\alpha_{6,2}=& {1 \over 648}\,(n^2-1)^2\,(m^2-1)^2\,n^2\,m^2, \cr
\tilde\alpha_{6,3}=& {1 \over 2160} \,(n^2-1)^2\,(m^2-1)^2
     ( 9\,{n^2}\,{m^2} - {n^2} - {m^2} -1), \cr
\tilde\alpha_{6,4}=& {1 \over 2160} \,(n^2-1)^2\,(m^2-1)^2 \, (n^2+1) \,
(m^2+1), \cr
\tilde\alpha_{6,5}=& {1 \over 75600} \, (n^2-1) \, (m^2-1) \,
\big(  516\,{n^4}\,{m^4}- 289\,({n^2}\,{m^4}+{n^4}\,{m^2}) \cr &
\,\,\,\,\,\,\,\,\,\,\,\,\,\,\,\,\,\,\,\,\,\,\,\,\,\,\,\,\,\,\,\,
\,\,\,\,\,\,\,\,\,\,\,\,\,\,\,\,\,\,\,\,\,\,\,\,
- 44\,{n^2}\,{m^2} + 5\,({n^4}+m^4) + 5\,({n^2}+m^2) + 5 \big), \cr
\tilde\alpha_{6,6}=& {1 \over 90720} \, (n^2-1) \, (m^2-1) \,
\big(53\,{n^4}\,{m^4}- 101\,({n^2}\,{m^4}+{n^4}\,{m^2}) \cr &
\,\,\,\,\,\,\,\,\,\,\,\,\,\,\,\,\,\,\,\,\,\,\,\,\,\,\,\,\,\,\,\,
\,\,\,\,\,\,\,\,\,\,\,\,\,\,\,\,\,\,\,\,\,\,\,\,
- 115\,{n^2}\,{m^2}- 24\,({n^4}+m^4) - 24\,({n^2}+m^2)  -24  \big), \cr
\tilde\alpha_{6,7}=& {1 \over 226800} \, (n^2-1) \, (m^2-1) \,
\big( 419\,{n^4}\,{m^4}+ 209\,({n^2}\,{m^4}+{n^4}\,{m^2}) \cr &
\,\,\,\,\,\,\,\,\,\,\,\,\,\,\,\,\,\,\,\,\,\,\,\,\,\,\,\,\,\,\,\,
\,\,\,\,\,\,\,\,\,\,\,\,\,\,\,\,\,\,\,\,\,\,\,\,
- {n^2}\,{m^2} + 20\,({n^4}+m^4)  + 20\,({n^2}+m^2) +20 \big), \cr
\tilde\alpha_{6,8}=& {1 \over 453600} \, (n^2-1) \, (m^2-1) \,
\big( 13\,{n^4}\,{m^4} + 13\,({n^2}\,{m^4}+{n^4}\,{m^2}) \cr &
\,\,\,\,\,\,\,\,\,\,\,\,\,\,\,\,\,\,\,\,\,\,\,\,\,\,\,\,\,\,\,\,
\,\,\,\,\,\,\,\,\,\,\,\,\,\,\,\,\,\,\,\,\,\,\,\,
+ 13\,{n^2}\,{m^2}  - 50\,({n^4}+m^4) - 50\,({n^2}+m^2) -50 \big), \cr
\tilde\alpha_{6,9}=& {1 \over 151200} \, (n^2-1) \, (m^2-1) \,
\big( 31\,{n^4}\,{m^4} + 31\,({n^2}\,{m^4} + {n^4}\,{m^2}) \cr &
\,\,\,\,\,\,\,\,\,\,\,\,\,\,\,\,\,\,\,\,\,\,\,\,\,\,\,\,\,\,\,\,
\,\,\,\,\,\,\,\,\,\,\,\,\,\,\,\,\,\,\,\,\,\,\,\,
+ 31\,{n^2}\,{m^2} + 10\,({n^4}+m^4) + 10\,({n^2}+m^2) + 10 \big).
\cr}
\eqn\alfatilde
$$
One can verify that the Vassiliev invariants which do not correspond to
primitive ones satisfy the relations [\numbers]:
$$
\eqalign{
\tilde\alpha_{4,1}=&\half \tilde\alpha_{2,1}^2, \cr
\tilde\alpha_{5,1}=&\tilde\alpha_{2,1}\tilde\alpha_{3,1}, \cr
\tilde\alpha_{6,1}=&{1\over 6}\tilde\alpha_{2,1}^3, \cr}
\qquad\qquad
\eqalign{
\tilde\alpha_{6,2}=&\half \tilde\alpha_{3,1}^2, \cr
\tilde\alpha_{6,3}=&\tilde\alpha_{2,1}\tilde\alpha_{4,2}, \cr
\tilde\alpha_{6,4}=&\tilde\alpha_{2,1}\tilde\alpha_{4,3}. \cr}
\eqn\uva
$$
As explained in [\numbers] these relations are a consequence of the
factorization property \factorization.

The quantities $\alpha_{ij}$, which are the ones which have a more direct
integral interpretation are easily obtained from the $\tilde\alpha_{ij}$ in
\alfatilde\ and the expansion for the unknot up to order six \unknot.
Since \unknot\ only contains even powers of $x$, for the invariants of
odd order one has:
$$
\eqalign{
\alpha_{3,1}=\tilde\alpha_{3,1}, \cr
\alpha_{5,2}=\tilde\alpha_{5,2}, \cr}
\qquad\qquad
\eqalign{
\alpha_{5,3}=\tilde\alpha_{5,3}, \cr
\alpha_{5,4}=\tilde\alpha_{5,4}. \cr}
\eqn\saturno
$$
As shown in [\numbers], the relations \uva\ for the compound invariants also
holds in this case:
$$
\eqalign{
\alpha_{4,1}=&\half \alpha_{2,1}^2, \cr
\alpha_{5,1}=&\alpha_{2,1}\alpha_{3,1}, \cr
\alpha_{6,1}=&{1\over 6}\alpha_{2,1}^3, \cr}
\qquad\qquad
\eqalign{
\alpha_{6,2}=&\half \alpha_{3,1}^2, \cr
\alpha_{6,3}=&\alpha_{2,1}\alpha_{4,2}, \cr
\alpha_{6,4}=&\alpha_{2,1}\alpha_{4,3}. \cr}
\eqn\uvados
$$
For the primitive ones of even order one finds:
$$
\eqalign{
\alpha_{2,1}=&{1\over 6}\,(n^2\, m^2-n^2-m^2), \cr
\alpha_{4,2}=&{1 \over 360} \,(9\,n^4\,
m^4-10\,(n^2\,m^4+n^4\,m^2)+(n^4+m^4) +10 n^2\,m^2), \cr
\alpha_{4,3}=& {1 \over 360} \,(n^4\, m^4-n^4-m^4), \cr}
$$
$$
\eqalign{
\alpha_{6,5}=& {1 \over 75600} \,
\big( 516\,{n^6}\,{m^6}-805\,(n^4\,m^6+n^6\,m^4)
+1050\,n^4\,m^4  \cr
&\,\,\,\,\,\,\,\,\,\,\,\,\,\,\,\,\,\,\,\,\,\,
+294\,(n^2\,m^6+n^6\,m^2)-245\,(n^2\,m^4+n^4\,m^2)\cr
&\,\,\,\,\,\,\,\,\,\,\,\,\,\,\,\,\,\,\,\,\,\,
-5\,(n^6+m^6)-49\,n^2\,m^2\big), \cr
\alpha_{6,6}=& {1 \over 90720} \,
\big( 53\,{n^6}\,{m^6}-154\,(n^4\,m^6+n^6\,m^4)
+140\,n^4\,m^4  \cr
&\,\,\,\,\,\,\,\,\,\,\,\,\,\,\,\,\,\,\,\,\,\,
+77\,(n^2\,m^6+n^6\,m^2)+14\,(n^2\,m^4+n^4\,m^2)\cr
&\,\,\,\,\,\,\,\,\,\,\,\,\,\,\,\,\,\,\,\,\,\,
+24\,(n^6+m^6)-91\,n^2\,m^2\big), \cr
\alpha_{6,7}=& {1 \over 226800} \,
\big( 419\,{n^6}\,{m^6}-210\,(n^4\,m^6+n^6\,m^4)
-189\,(n^2\,m^6+n^6\,m^2)\cr
&\,\,\,\,\,\,\,\,\,\,\,\,\,\,\,\,\,\,\,\,\,\,
+210\,(n^2\,m^4+n^4\,m^2)-20\,(n^6+m^6)-21\,n^2\,m^2\big),\cr
\alpha_{6,8}=& {1 \over 453600} \,
\big( 13\,{n^6}\,{m^6}-63\,(n^2\,m^6+n^6\,m^2)
+50\,(n^6+m^6)+63\,n^2\,m^2\big),\cr
\alpha_{6,9}=& {1 \over 151200} \,
\big( 31\,{n^6}\,{m^6}-21\,(n^2\,m^6+n^6\,m^2)
-10\,(n^6+m^6)+21\,n^2\,m^2\big).\cr }
\eqn\alfa
$$

The simple expressions found for the quantities $\alpha_{ij}$ are in clear
contrast with the complexity implicit in their integral form.  As it was
indicated above, the
$\alpha_{ij}$ constitute the form of the Vassiliev invariants which are more
simply related to their integral form as dictated from the Feynmann rules of
Chern-Simons gauge theory. These are known only for the first two invariants.
At order two one has [\gmm,\drortesis]:
$$
\eqalign{&
\alpha_{2,1}={1\over 4}\oint dx^{\mu_1}_1 \int^{x_1} dx^{\mu_2}_2
\int^{x_2} dx^{\mu_3}_3 \int^{x_3} dx^{\mu_4}_4
\Delta_{\mu_1\mu_3}(x_1-x_3) \Delta_{\mu_2\mu_4}(x_2-x_4) \cr
- &{1\over 16} \oint dx^{\mu_1}_1 \int^{x_1} dx^{\mu_2}_2
\int^{x_2} dx^{\mu_3}_3
\int_\erretres d^3 y \big( \Delta_{\mu_1\nu_1}(x_1-y)
\Delta_{\mu_2\nu_2}(x_2-y) \cr &
\,\,\,\,\,\,\,\,\,\,\,\,\,\,\,\,\,\,\,
\,\,\,\,\,\,\,\,\,\,\,\,\,\,\,\,\,\,\,
\,\,\,\,\,\,\,\,\,\,\,\,\,\,\,\,\,\,\,
\,\,\,\,\,\,\,\,\,\,\,\,\,\,\,\,\,\,\,
\,\,\,\,\,\,\,\,\,\,\,\,\,\,\,\,\,\,\,
\,\,\,\,\,\,\,\,\,\,\,\,\,\,\,\,\,\,\,
\times \Delta_{\mu_3\nu_3}(x_3-y)
\epsilon^{\nu_1\nu_2\nu_3}\big),\cr} \eqn\alfadosunoint
$$
while at order three [\numbers]:
$$
\eqalign{\alpha_{3,1}={1\over 8} &
\oint \,{dx_1^{\mu_1}}
\int^{x_1} \,{dx_2^{\mu_2}}
\int^{x_2} \,{dx_3^{\mu_3}}
\int^{x_3} \,{dx_4^{\mu_4}}
\int^{x_4} \,{dx_5^{\mu_5}}
\int^{x_5} \,{dx_6^{\mu_6}} \cr
 \big[ &
\Delta_{\mu_1\mu_4}(x_1-x_4)\Delta_{\mu_2\mu_6}(x_2-x_6)
\Delta_{\mu_3\mu_5}(x_3-x_5) \cr & +
\Delta_{\mu_1\mu_3}(x_1-x_3)\Delta_{\mu_2\mu_5}(x_2-x_5)
\Delta_{\mu_4\mu_6}(x_4-x_6)  \cr
&+\Delta_{\mu_1\mu_5}(x_1-x_5)\Delta_{\mu_2\mu_4}(x_2-x_4)
\Delta_{\mu_3\mu_6}(x_3-x_6)
\big] \cr
+{1\over 4}&
\oint \,{dx_1^{\mu_1}}
\int^{x_1} \,{dx_2^{\mu_2}}
\int^{x_2} \,{dx_3^{\mu_3}}
\int^{x_3} \,{dx_4^{\mu_4}}
\int^{x_4} \,{dx_5^{\mu_5}}
\int^{x_5} \,{dx_6^{\mu_6}} \cr
 \big[ &
\Delta_{\mu_1\mu_4}(x_1-x_4)\Delta_{\mu_2\mu_5}(x_2-x_5)
\Delta_{\mu_3\mu_6}(x_3-x_6)
\big] \cr
-{1\over 32}&
\oint \,{dx_1^{\mu_1}}
\int^{x_1} \,{dx_2^{\mu_2}}
\int^{x_2} \,{dx_3^{\mu_3}}
\int^{x_3} \,{dx_4^{\mu_4}}
\int^{x_4} \,{dx_5^{\mu_5}}
\int_{{\bf R}^3} d^3y \cr
 \big[ &
\Delta_{\mu_1\nu_1}(x_1-y)\Delta_{\mu_2\mu_5}(x_2-x_5)
\Delta_{\mu_3\nu_3}(x_3-y)\Delta_{\mu_4\nu_4}(x_4-y)
\epsilon^{\nu_1\nu_3\nu_4} \cr
&+ \Delta_{\mu_1\mu_3}(x_1-x_3)\Delta_{\mu_2\nu_2}(x_2-y)
\Delta_{\mu_4\nu_4}(x_4-y)\Delta_{\mu_5\nu_5}(x_5-y)
\epsilon^{\nu_2\nu_4\nu_5} \cr
&+\Delta_{\mu_1\nu_1}(x_1-y)\Delta_{\mu_2\mu_4}(x_2-x_4)
\Delta_{\mu_3\nu_3}(x_3-y)\Delta_{\mu_5\nu_5}(x_5-y)
\epsilon^{\nu_1\nu_3\nu_5} \cr
&+ \Delta_{\mu_1\nu_1}(x_1-y)\Delta_{\mu_2\nu_2}(x_2-y)
\Delta_{\mu_3\mu_5}(x_3-x_5)\Delta_{\mu_4\nu_5}(x_4-y)
\epsilon^{\nu_1\nu_2\nu_4} \cr
&+\Delta_{\mu_1\mu_4}(x_1-x_4)\Delta_{\mu_2\nu_2}(x_2-y)
\Delta_{\mu_3\nu_3}(x_3-y)\Delta_{\mu_5\nu_5}(x_5-y)
\epsilon^{\nu_2\nu_3\nu_5} \big] \cr
+{1\over 128}&
\oint \,{dx_1^{\mu_1}}
\int^{x_1} \,{dx_2^{\mu_2}}
\int^{x_2} \,{dx_3^{\mu_3}}
\int^{x_3} \,{dx_4^{\mu_4}}
\int_{{\bf R}^3} d^3y \int_{{\bf R}^3} d^3z\cr
\big[ &
\big( \Delta_{\mu_1\nu_1}(x_1-y)\Delta_{\mu_2\sigma_2}(x_2-z)
\Delta_{\mu_3\sigma_3}(x_3-z)\Delta_{\mu_4\nu_4}(x_4-y)
\cr &
\,\,\,\,\,\,\,\,\,\,\,\,\,\,\,\,\,\,\,
\,\,\,\,\,\,\,\,\,\,\,\,\,\,\,\,\,\,\,
\,\,\,\,\,\,\,\,\,\,\,\,\,\,\,\,\,\,\,
\,\,\,\,\,\,\,\,\,\,\,\,\,\,\,\,\,\,\,
\,\,\,\,\,\,\,\,\,\,\,\,\,\,\,\,\,\,\,
\,\,\,\,\,\,\,\,\,\,\,\,\,\,\,\,\,\,\,
\times\Delta_{\alpha\beta}(y-z)\epsilon^{\nu_1\alpha\nu_4}
\epsilon^{\sigma_3\beta\sigma_2} \big) \cr
&+\big( \Delta_{\mu_1\nu_1}(x_1-y)\Delta_{\mu_2\nu_2}(x_2-y)
\Delta_{\mu_3\sigma_3}(x_3-z)
\cr &
\,\,\,\,\,\,\,\,\,\,\,\,\,\,\,\,\,\,\,
\,\,\,\,\,\,\,\,\,\,\,\,\,\,\,\,\,\,\,
\,\,\,\,\,\,\,\,\,\,\,\,\,\,\,\,\,\,\,
\,\,\,\,\,\,\,\,\,\,\,\,\,\,\,\,\,\,\,
\times\Delta_{\mu_4\sigma_4}(x_4-z)
\Delta_{\alpha\beta}(y-z)\epsilon^{\nu_2\alpha\nu_1}
\epsilon^{\sigma_4\beta\sigma_3}\big) \big], \cr}
\eqn\alfatresuno
$$
where,
$$
\Delta_{\mu\nu}(x-y)= {1\over \pi} \epsilon_{\mu\sigma\nu}
{ (x-y)^\sigma \over |x-y|^3 }.
$$

We proposed in [\numbers] a normalization of the Vassiliev invariants which
seemed to be special in the sense that these invariants become integer-valued.
This normalization is fixed assuming that the Vassiliev invariants
$\tilde\alpha_{ij}$ for the trefoil knot do not vanish and then using those
non-vanishing values to normalize for the rest of the non-trivial knots
inserting appropriate factors. The resulting primitive invariants for a knot
$K$, which are labelled as $\beta_{ij}(K)$, are defined up to order six as:
$$
\eqalign{
\beta_{2,1}(K) &=    {\tilde\alpha_{2,1}(K) \over
\tilde\alpha_{2,1}(3_1)},\cr
\beta_{3,1}(K) &=    {\tilde\alpha_{3,1}(K) \over
\tilde\alpha_{3,1}(3_1)},\cr
\beta_{4,2}(K) &=  31{\tilde\alpha_{4,2}(K) \over
\tilde\alpha_{4,2}(3_1)},\cr
\beta_{4,3}(K) &=   5{\tilde\alpha_{4,3}(K) \over
\tilde\alpha_{4,3}(3_1)},\cr
\beta_{5,2}(K) &=  11{\tilde\alpha_{5,2}(K) \over
\tilde\alpha_{5,2}(3_1)},\cr
\beta_{5,3}(K) &=    {\tilde\alpha_{5,3}(K) \over
\tilde\alpha_{5,3}(3_1)},\cr}
\qquad\qquad
\eqalign{
\beta_{5,4}(K) &=    {\tilde\alpha_{5,4}(K) \over
\tilde\alpha_{5,4}(3_1)},\cr
\beta_{6,5}(K) &=5071{\tilde\alpha_{6,5}(K) \over
\tilde\alpha_{6,5}(3_1)},\cr
\beta_{6,6}(K) &=  29{\tilde\alpha_{6,6}(K) \over
\tilde\alpha_{6,6}(3_1)},\cr
\beta_{6,7}(K) &=1531{\tilde\alpha_{6,7}(K) \over
\tilde\alpha_{6,7}(3_1)},\cr
\beta_{6,8}(K) &=  17{\tilde\alpha_{6,8}(K) \over
\tilde\alpha_{6,8}(3_1)},\cr
\beta_{6,9}(K) &= 271{\tilde\alpha_{6,9}(K) \over
\tilde\alpha_{6,9}(3_1)},\cr}
\eqn\trebolbetas
$$
where $\tilde\alpha_{ij}(3_1)$ refers to the invariants corresponding
to the trefoil knot. The compound Vassiliev
invariants are then simply defined in such a way that they are products of
primitive ones: $$
\eqalign{
\beta_{4,1}=&\beta_{2,1}^2, \cr
\beta_{5,1}=&\beta_{2,1}\beta_{3,1}, \cr
\beta_{6,1}=&\beta_{2,1}^3, \cr}
\qquad\qquad
\eqalign{
\beta_{6,2}=&\beta_{3,1}^2, \cr
\beta_{6,3}=&\beta_{2,1}\beta_{4,2}, \cr
\beta_{6,4}=&\beta_{2,1}\beta_{4,3}. \cr}
\eqn\uvadosdos
$$

We will now compute the primitive Vassiliev invariants for torus knots
invariants in such a normalization. From \alfatilde\ follows:
$$
\eqalign{
\beta_{2,1}=& {1 \over 24}  \,( {n^2}-1 ) \,({m^2}-1),  \cr
\beta_{3,1}=& {1 \over 144} \,n\,m\,(  {n^2}-1 ) \,(  {m^2}-1 ), \cr
\beta_{4,2}=& {1 \over 240} \,( {n^2}-1 ) \, (m^2-1) \,
     ( 9\,{n^2}\,{m^2}-n^2-m^2-1 ), \cr
\beta_{4,3}=& {1 \over 240} \,(  {n^4}-1 ) \,(  {m^4}-1 ), \cr}
$$
$$
\eqalign{
\beta_{5,2}=& {1 \over 28800} \,n\,m \,( n^2-1 )\, ( m^2-1) \,
     ( 69\,{n^2}\,{m^2}  - 21\,\big({n^2} + {m^2})-11\big), \cr
\beta_{5,3}=& {1 \over 57600} \,n\,m \,( n^2-1 )\, ( m^2-1) \,
     ( 11\,{n^2}\,{m^2} + {n^2} + {m^2}-9), \cr
\beta_{5,4}=& {1 \over 7200} \,n\,m\,(  {n^4}-1 ) \,( {m^4}-1 ), \cr}
$$
$$
\eqalign{
\beta_{6,5}=& {1 \over 2520} \, (n^2-1) \, (m^2-1) \,
\big(  516\,{n^4}\,{m^4}- 289\,({n^2}\,{m^4}+{n^4}\,{m^2}) \cr &
\,\,\,\,\,\,\,\,\,\,\,\,\,\,\,\,\,\,\,\,\,\,\,\,\,\,\,\,\,\,\,\,
\,\,\,\,\,\,\,\,\,\,\,\,\,\,\,\,\,\,\,\,\,\,\,\,
- 44\,{n^2}\,{m^2} + 5\,({n^4}+m^4) + 5\,({n^2}+m^2) + 5 \big), \cr
\beta_{6,6}=& {1 \over 12096} \, (n^2-1) \, (m^2-1) \,
\big(53\,{n^4}\,{m^4}- 101\,({n^2}\,{m^4}+{n^4}\,{m^2}) \cr &
\,\,\,\,\,\,\,\,\,\,\,\,\,\,\,\,\,\,\,\,\,\,\,\,\,\,\,\,\,\,\,\,
\,\,\,\,\,\,\,\,\,\,\,\,\,\,\,\,\,\,\,\,\,\,\,\,
- 115\,{n^2}\,{m^2}- 24\,({n^4}+m^4) - 24\,({n^2}+m^2)  -24  \big), \cr
\beta_{6,7}=& {1 \over 10080} \, (n^2-1) \, (m^2-1) \,
\big( 419\,{n^4}\,{m^4}+ 209\,({n^2}\,{m^4}+{n^4}\,{m^2}) \cr &
\,\,\,\,\,\,\,\,\,\,\,\,\,\,\,\,\,\,\,\,\,\,\,\,\,\,\,\,\,\,\,\,
\,\,\,\,\,\,\,\,\,\,\,\,\,\,\,\,\,\,\,\,\,\,\,\,
- {n^2}\,{m^2} + 20\,({n^4}+m^4)  + 20\,({n^2}+m^2) +20 \big), \cr
\beta_{6,8}=& {1 \over 25200} \, (n^2-1) \, (m^2-1) \,
\big( 13\,{n^4}\,{m^4} + 13\,({n^2}\,{m^4}+{n^4}\,{m^2}) \cr &
\,\,\,\,\,\,\,\,\,\,\,\,\,\,\,\,\,\,\,\,\,\,\,\,\,\,\,\,\,\,\,\,
\,\,\,\,\,\,\,\,\,\,\,\,\,\,\,\,\,\,\,\,\,\,\,\,
+ 13\,{n^2}\,{m^2}  - 50\,({n^4}+m^4) - 50\,({n^2}+m^2) -50 \big), \cr
\beta_{6,9}=& {1 \over 5040} \, (n^2-1) \, (m^2-1) \,
\big( 31\,{n^4}\,{m^4} + 31\,({n^2}\,{m^4} + {n^4}\,{m^2}) \cr &
\,\,\,\,\,\,\,\,\,\,\,\,\,\,\,\,\,\,\,\,\,\,\,\,\,\,\,\,\,\,\,\,
\,\,\,\,\,\,\,\,\,\,\,\,\,\,\,\,\,\,\,\,\,\,\,\,
+ 31\,{n^2}\,{m^2} + 10\,({n^4}+m^4) + 10\,({n^2}+m^2) + 10 \big).
\cr}
\eqn\betas
$$
In the following section we will describe some of the properties of these
invariants.

\chapter{Properties of the Vassiliev Invariants for Torus Knots}

Given the Vassiliev invariants \betas\ the first question that one would like
to answer is how many of them are needed to distinguish torus knots. It was
shown in [\trapp] that for torus knots of the form $(2,2p+1)$ the Vassiliev
invariants of second and third orders are enough. We will prove now that this
result holds for arbitrary torus knots. In other words, we will prove that if
two torus knots $\{n,m\}$ (with $(n,m)=1$) and $(n',m')$ (with $(n',m')=1$)
have the same first two Vassiliev invariants,
$$
\eqalign{
\beta_{2,1}(n,m) = \beta_{2,1}(n',m'), \cr
\beta_{3,1}(n,m) = \beta_{3,1}(n',m'), \cr}
\eqn\marte
$$
then $(n',m')=(n,m)$, $(n',m')=(m,n)$, $(n',m')=(-n,-m)$, or
$(n',m')=(-m,-n)$. First notice that from the explicit form of
$\beta_{2,1}$ and $\beta_{3,1}$ in \betas, one finds from \marte:
$$
m\,n = m'\,n',   \,\,\,\,\,\,\,\,\,\,\,\,\,\,\,\,\,
{\hbox{\rm and}}  \,\,\,\,\,\,\,\,\,\,\,\,\,\,\,\,\,
n^2+m^2 = {n'}^2+{m'}^2.
\eqn\venus
$$
Let $\gamma$ be a rational number such that  $m'=\gamma\,n$. From
the first relation in \venus\ follows that also $m=\gamma\,n'$. Then, from the
second relation in \venus\ one gets:
$$
n^2(\gamma^2-1) = {n'}^2(\gamma^2-1),
\eqn\tierra
$$
which implies that either $n=\pm n'$ and $\gamma$ arbitrary, or
$\gamma=\pm 1$ and $n$
and $n'$ any. In the first of these two cases, after considering \venus, one
ends with either $n=n'$ and $m=m'$, or $n=-n'$ and $m=-m'$. In the second
case, similarly, one concludes that either $m=n'$ and $n=m'$, or $m=-n'$ and
$n=-m'$. The four possibilities correspond to the same torus knot and
therefore we have proved that the first two Vassiliev invariants distinguish
all torus knots.

This analysis shows that the first two Vassiliev invariants provide a
one-to-one map between the set ${\cal T}=\{ n,m \in {\bf
Z}\;|\;(n,m)=1,\;n>|m|\}$
 and its image  through the function $(\beta_{2,1}(n,m),\beta_{3,1}(n,m))$.
{}From the form of $\beta_{2,1}$ and $\beta_{3,1}$ in \betas\ one finds that
${\cal T}$ gets mapped around a curve in the $(\beta_{2,1},\beta_{3,1})$ plane
whose asymptotic form (large values of $\beta_{2,1}$ and $\beta_{3,1}$) is:
$$
\beta_{3,1}^2={2\over 3} \beta_{2,1}^3.
\eqn\curva
$$

It is important to remark that our result for $\beta_{3,1}$ agrees with
the Vassiliev invariant of order three for
torus knots of the form $(2,2p+1)$ presented in [\trapp]: $v_3=p^3-p$.
To verify this one must take into account that Vassiliev invariants are
defined up to global factors and addition of a linear relation of lower order
Vassiliev invariants. It turns out that:
$$
v_3 = 3 (\beta_{3,1}(2,2p+1)-\beta_{2,1}(2,2p+1) ) =p^3-p.
\eqn\jupiter
$$
Also, it would be interesting to study if for the lower order cases (two
and three) the methods developped in
\REF\sst{A.C. Hirshfeld and U. Sassenberg, ``Derivation of the total twist
from Chern-Simons theory", Preprint DO-TH 95/02, hep-th/9502088}
[\sst] lead to  results equivalent to ours.

Another important issue that we will address here is to work out which of the
Vassiliev invariants which we have obtained contain  the topological
information. In other words, we will study at each order how many
independent invariants are taking into account that linear combinations of
lower-order invariants can be added.
Certainly, at orders two and three there is a
single one for each. At order four one finds that,
$$
\beta_{4,2} = 4 \beta_{4,3} + 12 \beta_{2,1}^2 - \beta_{2,1},
\eqn\pluton
$$
and therefore there is only one independent invariant. At order five one
obtains,
$$
\eqalign{
\beta_{5,2} =& 6 \beta_{5,3} + {27\over 5} \beta_{2,1}\beta_{3,1}
 - {2\over 5} \beta_{3,1}, \cr
\beta_{5,3} =& {3\over 4} \beta_{5,3} + {3\over 10} \beta_{2,1}\beta_{3,1}
 - {1\over 20} \beta_{3,1}, \cr}
\eqn\moon
$$
and thus, again, there is only one independent invariant. At order six,
however one finds that the number of independent invariants is two. Indeed,
from the relations,
$$
\eqalign{
\beta_{6,5} =& {58\over 9} \beta_{6,9} - {80\over 3} \beta_{4,3}
+{41\over 9}\beta_{2,1}-{680\over 3}\beta_{2,1}\beta_{4,3}
+5280\beta_{3,1}^2-{2080\over 3}\beta_{2,1}^3, \cr
\beta_{6,6} =& -{5\over 12} \beta_{6,9} - {5\over 3} \beta_{4,3}
+{1\over 4}\beta_{2,1}-{10}\beta_{2,1}\beta_{4,3}
+240\beta_{3,1}^2-{40}\beta_{2,1}^3, \cr
\beta_{6,7} =& {9\over 2} \beta_{6,9} - {5} \beta_{4,3}
+{1\over 2}\beta_{2,1}
+432\beta_{3,1}^2-{96}\beta_{2,1}^3, \cr}
\eqn\moon
$$
and the fact that $\beta_{6,8}$ can not be written as a linear relation of
the $\beta_{ij}$ which appear on the right hand side of \moon\ one concludes
that $\beta_{6,8}$ and $\beta_{6,9}$ can be taken as the independent
invariants at order six.  Notice that up to the order which we have studied
the independent Vassiliev
invariants for torus knots can be chosen for $i>3$ as the ones associated to
the Casimirs shown in Fig. 2. As discussed in [\numbers], these Casimirs are
the building blocks of the group factors shown in Fig. 1. One would like to
know if the correspondence found up to order six holds in general.

The result obtained for $\beta_{2,1}$ allows to conclude that many torus
knots are not Lissajous knots. It has been shown recently
\REF\lisa{M.G.V. Bogle, J.E. Hearst, V.F.R. Jones and
L. Stoilov\journal\knot&3(94)121} [\lisa]
that the Arf, Kervaire, or Robertello invariant for a Lissajous knot is zero.
On the other hand, it is also known [\gmm] that this invariant is just
$\beta_{2,1}$ mod 2. Thus, if $\beta_{2,1}$ is odd for a torus knot $\{n,m\}$
one can state that it is not a Lissajous knot. From the formula for
$\beta_{2,1}$ in \betas\ one can verify that indeed the trefoil knot ($n=2$,
$m=3$) is not a Lissajous knot. Our result for the Vassiliev invariant
$\beta_{2,1}$ in \betas\ allows to conclude in general that if
$(n^2-1)(m^2-1)/24$ is odd the torus knot $\{n,m\}$ is not a Lissajous knot.

One of the most interesting property of the Vassiliev
invariants $\beta_{ij}$ is
that they seem to be integer-valued. It was observed in [\numbers] that there
seemed to be a special normalization for the Vassiliev invariants such that
they become integer-valued. One could say that such observation was not
well-funded because after all in [\numbers] only a finite set of knots were
considered. However, in this paper we have studied an infinite subset of
knots and one seems to find the same property. It is proved in the appendix
that the $\beta_{ij}$ in \betas, which are polynomials in $n$ and $m$, are
integer-valued when $n$ and $m$ are coprime integers, $(n,m)=1$, up to order
four.  For orders
$i=5,6$ we have numerical evidence that this feature also hold but we do
not have yet a proof.
These polynomials share the property that there exist at least one
value (then there are infinitely many) of the pair $n,\,m$, such that
$(n,m)\neq 1$, for which $\beta_{ij}$ is not integer-valued.

The invariants $\beta_{ij}$ regarded as polynomials in $n$ and $m$ are very
interesting by themselves from the point of view of number theory. These are
polynomials which are integer-valued when $(n,m)=1$ but  fail for some  $n,\,m$
when $(n,m)\neq 1$. Given a polynomial with these features one
can always construct a new one adding a new polynomial at most of the same
degree with integer coefficients. Therefore, they are defined modulo
the ring of polynomials in $n$ and $m$ with integer coefficients. The prototype
for these polynomials is actually the one leading to the Gordian number of a
torus knot. It has the form $(n-1)(m-1)/2$ (for $n,m >0$) and according to
Milnor's conjecture
\REF\milnor{J. Milnor, ``Singular Points of Complex Hypersurfaces",
Princeton University Press, 1968} [\milnor]
corresponds to the uncrossing number of a torus knot
$\{n,m\}$.  Among the set of symmetric polynomials in $n$ and $m$ which only
vanish for $n=1$ and $m=1$, the polynomial $(n-1)(m-1)/2$ is, up to a global
sign, the only one which is integer-valued when $(n,m)=1$ but  it is not
integer-valued  for some $n,\,m$ such that $(n,m)\neq 1$, modulo the ring of
polynomials  in $n$ and $m$ with integer coefficients.

Relative to the Gordian number, notice that our results show that
$(n-1)(m-1)/2$ is not a Vassiliev invariant. If Milnor's
conjecture holds this is consistent with the fact proven in [\dean,\trapp]
that the uncrossing number is not an invariant of finite type.

According to \betas, the
polynomials which seem to be relevant are symmetric polynomials in $n$ and
$m$ which vanish for $n=\pm 1$ and $m=\pm 1$, and have the property described
above: they are integer-valued when $(n,m)=1$, and they fail to be an integer
for some $n,\,m$ such that $(n,m)\neq 1$. Actually, the set of polynomials in
which one is interested  is more restricted, for even orders they are
invariant under  $n\rightarrow -n$ (or $m\rightarrow -m$), while for odd
orders they change sign under such a transformation. The properties and
structure of this set of polynomials deserve to be studied
from a general point of view.

We would like to end pointing out that one of the important consequences of
knowing the explicit form of the Vassiliev invariants $\beta_{ij}$ at least
for torus knots is that it might shed light towards
their  general interpretation.

\vskip0.5cm

\ack

We would like to thank A. V. Ramallo and C. Lozano for very helpful
discussions. This work was supported in part by DGICYT under grant
PB93-0344 and  by CICYT under grant AEN94-0928.

\appendix

In this appendix we will prove that the Vassiliev invariants for torus knots
as given in \betas\ are integer-valued for orders $i=2,3,4$.

\newdimen\tableauside\tableauside=1.0ex
\newdimen\tableaurule\tableaurule=0.4pt
\newdimen\tableaustep
\def\phantomhrule#1{\hbox{\vbox to0pt{\hrule height\tableaurule width#1\vss}}}
\def\phantomvrule#1{\vbox{\hbox to0pt{\vrule width\tableaurule height#1\hss}}}
\def\sqr{\vbox{%
  \phantomhrule\tableaustep
  \hbox{\phantomvrule\tableaustep\kern\tableaustep\phantomvrule\tableaustep}%
  \hbox{\vbox{\phantomhrule\tableauside}\kern-\tableaurule}}}
\def\squares#1{\hbox{\count0=#1\noindent\loop\sqr
  \advance\count0 by-1 \ifnum\count0>0\repeat}}
\def\tableau#1{\vcenter{\offinterlineskip
  \tableaustep=\tableauside\advance\tableaustep by-\tableaurule
  \kern\normallineskip\hbox
    {\kern\normallineskip\vbox
      {\gettableau#1 0 }%
     \kern\normallineskip\kern\tableaurule}%
  \kern\normallineskip\kern\tableaurule}}
\def\gettableau#1 {\ifnum#1=0\let\next=\null\else
  \squares{#1}\let\next=\gettableau\fi\next}

\tableauside=1.2ex
\tableaurule=0.4pt

{\bf Proposition 1.} If $n$ and $m$ are two
coprime integers, $(n,m)=1$, then $\beta_{2,1}=(n^2-1)(m^2-1)/24$ is
integer-valued.

\noindent
{\it Proof:} Let $n$ and $m$ be two integers such that $(n,m)=1$. Being
coprime there is at least one of them which is odd. Let this be $n$. Since the
polynomial entering $\beta_{2,1}$ is symmetric in $n$ and $m$ this can be done
without loss of generality. If $n$ is odd one has that either $n-1=0$ mod
$2$ and
$n+1=0$ mod $4$, or $n+1=0$ mod $2$ and $n-1=0$ mod $4$. In either case, $n^2-1
= 0$ mod $8$. If in addition $n$ has the property of being $n=1,2$ mod $3$ one
has that either $n-1=0$ mod $3$, or $n+1=0$ mod $3$. In either case $n^2-1=0$
mod $3$ and therefore $n^2-1=0$ mod $24$, probing the proposition for this
case.
For the case which is left, $n=0$ mod $3$, one has that since $(n,m)=1$,
$m=1,2$
mod $3$ which implies as before $m^2-1=0$ mod $3$. Then, $(n^2-1)(m^2-1)=0$
mod 24, and the propostion is proven.  $\tableau{1}$

{\bf Proposition 2.} If $n$ and $m$ are two
coprime integers, $(n,m)=1$, then $\beta_{3,1}=nm(n^2-1)(m^2-1)/144$ is
integer-valued.

\noindent
{\it Proof:} Let $n$ and $m$ be two integers such that $(n,m)=1$. Being
coprime, as in the previous proof one can always choose $n$ to be odd without
loss of generality. As shown there one has then that $n^2-1 = 0$ mod
$8$ and, if in addition $n=1,2$ mod $3$, one has that $n^2-1=0$ mod
$24$. Let us consider the possible cases for $m$ under this situation.
We will consider the other situation (the one in which $n=0$ mod $3$ below).
If $m$ is odd and $m=1,2$ mod 3, the same arguments as before lead to
$m^2-1=0$ mod $24$, and therefore $(n^2-1)(m^2-1)=0$ mod $24^2$. If it where
the case in which $m=0$ mod $3$ one would have $m(m^2-1)=0$ mod 24, and
therefore $m(n^2-1)(m^2-1)=0$ mod $24^2$. Finally, if $m$ is even either
$m+1=0$ mod 3, or $m-1=0$ mod 3, being then $m(m^2-1)=0$ mod 6, and therefore
$m(n^2-1)(m^2-1)=0$ mod $144$. Let us consider now the second situation. If
$n=0$ mod 3 one has that $n(n^2-1)=0$ mod 24. If $m$ is odd, since $(n,m)=1$
one must have $m=1,2$ mod 3 and, as before, $m^2-1=0$ mod 24, obtaining then
$n(n^2-1)(m^2-1)=0$ mod $24^2$. If $m$ is even, again $m(m^2-1)=0$ mod 6, and
then $nm(n^2-1)(m^2-1)=0$ mod 144. This ends the proof of the proposition.
$\tableau{1}$

{\bf Proposition 3.} If $n$ and $m$ are two
coprime integers, $(n,m)=1$, then $\beta_{4,3}=(n^4-1)(m^4-1)/240$ is
integer-valued.

\noindent
{\it Proof:} Let $n$ and $m$ be two integers such that $(n,m)=1$. Being
coprime, as in the previous proof one can always choose $n$ to be odd without
loss of generality. As shown there one has then that $n^2-1 = 0$ mod
$8$. We will consider four possible situations:

\noindent
{\sl a)} $n=1,2$ mod 3 and $n=1,2,3,4$ mod 5. On the one hand, as in the
previous proposition one has that $n^2-1=0$ mod 24. On the other hand,
$n^2+1=0$
mod 2 and, if $n=2,3$ mod 5, $n^2+1=0$ mod 5, while if $n=1,4$ mod 5, $n^2-1=0$
mod 5. Then $(n^2-1)(n^2+1)=0$ mod 240.

\noindent
{\sl b)} $n=1,2$ mod 3 and $n=0$ mod 5. On the one hand, as before,
$n^2-1=0$ mod 24, and $n^2+1=0$ mod 2.
On the other hand, since $(n,m)=1$, $m=1,2,3,4$ mod 5. As in the previous
situation, then $(m^2-1)(m^2+1)=0$ mod 5. Therefore, one concludes that
$(n^2-1)(n^2+1)(m^2-1)(m^2+1)=0$ mod 240.

\noindent
{\sl c)} $n=0$ mod 3 and $n=1,2,3,4$ mod 5. In this case one has only that
$n^2-1=0$ mod 8 and $n^2+1=0$ mod 2, while, again, $(n^2-1)(n^2+1)=0$ mod 5.
As $(n,m)=1$ one must have $m=1,2$ mod 3, which implies that $m^2-1=0$ mod 3.
Then, $(n^2-1)(n^2+1)(m^2-1)(m^2+1)=0$ mod 240.

\noindent
{\sl d)} $n=0$ mod 3 and $n=0$ mod 5. In this case one has only that
$n^2-1=0$ mod 8 and $n^2+1=0$ mod 2. On the other hand, as
$(n,m)=1$ one must have $m=1,2$ mod 3, which implies that $m^2-1=0$ mod 3, and
$m=1,2,3,4$ mod 5, which implies that $(m^2-1)(m^2+1)=0$ mod 5.
Then, $(n^2-1)(n^2+1)(m^2-1)(m^2+1)=0$ mod 240. $\tableau{1}$

{\bf Proposition 4.} If $n$ and $m$ are two
coprime integers, $(n,m)=1$, then
$\beta_{4,2}=(n^2-1)(m^2-1)(9n^2m^2-n^2-m^2-1)/240$ is integer-valued.

\noindent
{\it Proof:} This proposition follows from
equation \pluton\ and Propositions 1 and 3. $\tableau{1}$

\refout

\end


\centerline{\bf VASSILIEV INVARIANTS FOR TORUS KNOTS}
\vskip1cm
\centerline{M. ALVAREZ}
\centerline{Center for Theoretical Physics}
\centerline{Massachusetts Institute of Technology}
\centerline{Cambridge, MA 02139, USA}
\centerline{{\it e-mail:} {\tt alvarez@mit} }
\vskip0.5cm
\centerline{and}
\vskip0.5cm
\centerline{J.M.F. LABASTIDA}
\centerline{Departamento de F\'\i sica de Part\'\i culas}
\centerline{Universidade de Santiago}
\centerline{E-15706 Santiago de Compostela, Spain}
\centerline{{\it e-mail:} {\tt labastida@gaes.usc.es} }

\vskip2cm

\centerline{ABSTRACT}

{\narrower\smallskip\ninepoint\baselineskip10pt
Vassiliev invariants up to order six for arbitrary
torus knots $\{ n , m \}$, with $n$ and $m$ coprime integers are computed.
These invariants are polynomials in $n$ and $m$ whose degree coincide with
their order. Furthermore, they turn out to be integer-valued in a normalization
previously proposed by the authors.
\smallskip}

{\narrower\smallskip\ninepoint\baselineskip10pt\noindent
{\it Keywords:} torus, knot, Vassiliev, Chern-Simons.
\smallskip}

